\documentstyle[12pt]{article}
\newtheorem{lemma}{Lemma}
\newtheorem{theorem}{Theorem}

\newcommand{\QED}{{\hfill$\Box$\medskip}}

\newtheorem{corollary}{Corollary}

\begin{document}
\begin{center}
{\Large \bf The Yamabe invariant of simply connected manifolds}\\
\vspace{1 cm}
Jimmy Petean\\
Max-Planck Institut f\"ur Mathematik,\\
Bonn, Germany\\
\end{center}

\vspace{1 cm}
\begin{abstract}
Let $M$ be any  simply connected smooth compact manifold
of dimension $n\geq 5$. We prove that the Yamabe 
invariant of $M$ is non-negative.
This is equivalent to say that the 
infimum, over the space of all Riemannian metrics
on $M$, of the $L^{n/2}$ norm of the  scalar
curvature is zero. 
\end{abstract}

\section{Introduction}

A classical problem in differential geometry is the determination
of the family of functions that can be obtained  as the scalar curvature
of a Riemannian metric in some fixed smooth manifold $M$.
As one might expect, the features of the problem in low dimensions
are very different from the high-dimensional case.
The classical Uniformization Theorem assures that any Riemannian metric
on a 2-dimensional manifold is conformal to a metric of constant
scalar curvature (in dimension two, this is the same as constant sectional
curvature). For a compact surface, 
the Gauss-Bonnet  formula  then shows that the sign of such a  constant
is determined by the topology of the surface. 
Moreover, if we restrict to metrics of unit  volume, then the value
of the constant is completely determined by the topology (it is $4\pi \chi$,
where $\chi $ is the Euler characteristic of the surface).

In higher dimensions, we can still deform any metric to a metric of 
constant scalar
curvature. Namely, given
any Riemannian metric $g$ on a compact smooth manifold $M$, there is a metric
conformal to $g$ which has constant scalar curvature. This is the
well-known Yamabe problem. Yamabe first stated this result in \cite{Yamabe},
but his proof was not correct, as pointed out by Trudinger in \cite{Trudinger}.
The argument was  completed in several steps by Trudinger
\cite{Trudinger}, Aubin \cite{Aubin}
and Schoen \cite{Schoen2}. Note that  in dimensions greater than two the
obstruction to the existence of metrics with a determined sign is much
weaker. For instance, every compact manifold of dimension at least three admits
a metric of constant negative scalar curvature. 
There are 
well-known obstructions, though, to the existence of metrics of positive
or vanishing scalar curvature.

If a function is expressed as the scalar curvature of a metric on a manifold
$M$, any positive multiple of the function can be  obtained by rescaling the 
metric. In order to avoid these ``trivial'' variations, and therefore to get
a meaningful measure of the possible ``size'' of the scalar curvature function,
it is reasonable to restrict 
to metrics of unit volume.
This was first considered by O. Kobayashi in \cite{Kobayashi}, where he
introduced what we will call the {\it Yamabe invariant} of a compact smooth
manifold $M$. First consider a fixed conformal class of metrics 
${\cal C}$ on $M$,  and
let the  Yamabe constant of $(M, {\cal C})$  be

$$Y(M,{\cal C})=\inf\limits_{g\in {\cal C}} 
\frac{\int\limits_M s_g \ dvol_g}{(Vol_g (M))^{2/n}}.$$

The solution to the Yamabe problem is precisely achieved by showing that
this infimum is always attained by a smooth metric, which necessarily
has constant scalar curvature. The {\it Yamabe invariant} is then defined by

$$Y(M)=\sup\limits_{\cal C} Y(M,{\cal C}),$$

\noindent
where the supremum is taken 
over all conformal classes of metrics on $M$. Note that
this invariant is also frequently called the {\it sigma constant} of $M$
\cite{Schoen3}. 

Note that the invariant is readily computable in dimension two from the
Gauss-Bonnet formula. In dimension three, the computation of the invariant
(of manifolds for which the invariant is non-positive)
would follow from Anderson's program for the hyperbolization 
conjecture \cite{Anderson}. Computations of the invariant in dimension
four have been carried out  by  LeBrun in \cite{LeBrun2}; he computed the
invariants of all compact complex surfaces of K\"{a}hler type
which do not admit metrics of positive scalar curvature. See also
\cite{Gursky, LeBrun3, Petean2} 
for other computations of the invariant in dimension four.

We will be concerned in this paper with 
simply connected manifolds of dimension greater than four.
The study  of manifolds which admit metrics of 
positive scalar curvature has proved
to be very interesting and deep. See for instance the work by Gromov
and Lawson
\cite{Gromov}, Schoen
and Yau \cite{Schoen} and Stolz \cite{Stolz}. In this
last paper, it is completely determined which simply connected  compact
manifolds of dimension greater than four admit metrics of positive scalar 
curvature. 
Recall that the Yamabe invariant of $M$ is positive
if and only if $M$ admits a metric
of positive scalar curvature (see for instance \cite{Schoen3}).

We will prove:

\begin{theorem}: Every simply connected smooth compact manifold of dimension
greater than four has non-negative Yamabe invariant.
\end{theorem}

Note that for a manifold $M$ which does not admit positive scalar curvature
metrics, the Yamabe invariant can also be computed by the formula

$$Y(M)= - \inf\limits_{\cal M} 
{\left( \int\limits_M |s_g |^{n/2} dvol_g \right)}^{2/n}$$

\noindent
where ${\cal M}$ is the space of all Riemannian metrics on $M$
(see for instance \cite{LeBrun2}).

It is also well-known that  if a compact manifold 
(of dimension at least three) admits a metric of positive scalar
curvature then it also admits a scalar-flat metric. We can therefore
rephrase the previous theorem as:

\begin{corollary}: For every simply connected smooth compact 
manifold M of dimension greater
than four,

$$\inf\limits_{\cal M} \int\limits_M |s_g |^{n/2} dvol_g =0.$$

\end{corollary}

It follows from these results that dimension four is quite exceptional 
in terms of the Yamabe invariant. Note that Theorem 1 is obviously true
in dimension two and clearly expected to be  true in the  3-dimensional
case (either from the Poincare Conjecture  or from Anderson's program
for the Hyperbolization Conjecture
\cite{Anderson}). But it is definitely not true in
dimension four  (see the work of LeBrun in \cite{LeBrun, LeBrun2}); moreover, 
it seems likely to be the case, from LeBrun's computations, 
that  the  {\it generic} simply 
connected compact four-manifold has strictly negative Yamabe
invariant.

\section{The Spin Cobordism Ring}

The aim of this section is the proof of Theorem 2 below. We will need to recall
some basic facts about the spin cobordism ring. Details can be found in
\cite{Lawson, Stong}.
Of course, we will need to consider in this section
manifolds with boundary. And hence, in this section,  
we will call a manifold $X$ 
{\it closed} if it is compact and without boundary (in the other sections
we always assume that the manifolds have no boundary).

For a spin manifold we will mean a smooth oriented manifold with a fixed
spin structure on its tangent bundle.
Let $X$ be a spin manifold with boundary. The spin structure on $X$  
induces a canonical spin structure on the boundary of $X$ (see \cite{Milnor}),
and two closed  
spin manifolds $X_1$ and $X_2$ of dimension $n$ 
are called {\it spin cobordant} if
there is a compact spin manifold $X$ (of dimension $n+1$) so that 
$\partial X$ is, as a spin manifold, 
the disjoint union of $X_1$ and  $-X_2$ (the minus meaning  
that the orientation is reversed).
The set of equivalence classes of $n$-dimensional closed spin manifolds
under this relation is  called the $n$-dimensional  
{\it spin cobordism group}, 
${\Omega}_n^{Spin}$. It is an Abelian group, with the sum given by  
the disjoint union or equivalently by the connected sum (if $X$ and
$Y$ are connected spin manifolds, then $X\# Y$ inherits a spin structure
so that $X\# Y$ is spin cobordant to the disjoint union of $X$ and $Y$).

Throughout this section
we will denote by $[X]$ the class of the closed spin manifold $X$
in the spin cobordism group.

The product of manifolds gives a ring structure to ${\Omega}_*^{Spin}$.
This is called the {\it spin cobordism  ring}. Its structure has been 
determined 
by D.W. Anderson, E.H. Brown and F.P. Peterson 
in \cite{Anderson3}. It is proved there that two closed
spin manifolds of dimenion $n$ are spin cobordant if and only if they
have the same Stiefel-Whitney and KO characteristic numbers. A particular
example of these  characteristic numbers, which will play an important role
in this article, is the $\alpha$-homomorphism,

$$\alpha :  {\Omega}_*^{Spin} \rightarrow KO_*  (pt),$$

\noindent
first introduced by Atiyah and Milnor (see  \cite{Milnor}). It is important
for us that $\alpha$ is a ring homomorphism; i.e. for any connected closed spin
manifolds $X$ and $Y$, $\alpha [X\# Y]=\alpha [X]+\alpha [Y]$ and
$\alpha [X\times Y]=\alpha [X] .\alpha [Y]$.

Now recall that $KO_n (pt)$ vanishes for
$n$=3,5,6 and 7, while $KO_1 (pt)$ and $KO_2 (pt)$ are isomorphic to
${\bf Z}_2 $, and
$KO_4 (pt)$ and $KO_8 (pt)$ are isomorphic to ${\bf Z}$.
Recall also that multiplication by
a generator of $KO_8 (pt)$  gives an isomorphism between
$KO_n  (pt)$ and $KO_{n+8} (pt)$. Moreover, $\alpha $ is an extension of
the $\hat{A}$-genus in the sense that after  picking suitable generators
of $KO_{8k}$ and $KO_{8k+4}$, $\alpha$ is exactly the $\hat{A}$-genus  in
dimensions 8k and one half of the $\hat{A}$-genus  in dimensions 8k+4.

The homomorphism $\alpha$  plays  a central role in the study of the scalar
curvature  of compact spin manifolds.
Hitchin proved in \cite{Hitchin} (generalizing the now classical work
of Lichnerowicz \cite{Lichnerowicz})
that if the compact spin manifold $X$ admits a metric of positive
scalar curvature, then $\alpha [X]=0$. The converse of  this result 
is also true for simply connected manifolds. It
was conjectured by Gromov and Lawson 
in \cite{Gromov}, and proved in different cases by
Gromov and Lawson \cite{Gromov}, 
Miyazaki \cite{Miyazaki}, Rosenberg \cite{Rosenberg}
and finally (in general) by Stolz \cite{Stolz}. More precisely, Stolz
proved that any cobordism class in 
the kernel of  the homomorphism $\alpha$ can be represented by a connected 
spin manifold which admits metrics of positive scalar curvature (the
result then follows from the arguments  in \cite{Gromov}).
Similarly, to prove Theorem 1, we will need to show that every element
in the spin cobordism groups can be represented by a manifold with non-negative
Yamabe invariant. We will prove this statement for every dimension in 
Theorem 2, although we will only need the result for dimensions greater than 
four.
We will prove Theorem 2 now and then use it to prove Theorem  1
in the following section.

\vspace{.2cm}

We will need  to use the following obvious fact:

\begin{lemma}: Suppose that $M$ and $N$ are closed smooth manifolds and
that the Yamabe invariant of $M$ is non-negative. Then 
the Yamabe invariant of $M\times N$ is non-negative.
\end{lemma}

Proof: Given any $\epsilon >0$ we will show that there is a metric on
$M\times N$ of unit volume and constant scalar curvature bounded below by 
$-\epsilon$. This implies the lemma.

Every compact smooth manifold admits a metric of constant scalar curvature.
After rescaling, we can then pick a metric $g$  on $N$  of volume $V$ 
whose scalar curvature  is a constant greater than 
$-\epsilon /2$. The fact that the Yamabe invariant of $M$ is
non-negative  assures that there
is a metric $h$ on $M$ of volume 1/$V$ and whose scalar curvature 
is a constant greater than $-\epsilon /2$. 
The product metric $g+h$ on $N\times M$ has unit volume  
and scalar curvature bounded below by $-\epsilon$
as required.

\QED

We can now prove:

\begin{theorem}: Every element in the spin cobordism group 
${\Omega}_n^{Spin}$  can be represented by a connected spin manifold
with non-negative Yamabe invariant.
\end{theorem}

Proof: The statement is trivial in the case $n=1$. When $n=2$, we have that
${\Omega}_2^{Spin}\ $is isomorphic to 
${\bf Z}_2$, and its two elements are represented
respectively by a torus (the product of two copies of the 
spin structure on $S^1$ that is not a spin boundary), 
whose Yamabe invariant is zero,
and by the 2-sphere with its canonical
spin structure (the Yamabe invariant of the 2-sphere is 8$\pi$).
When $n$=3,5,6 or 7, the spin cobordism group ${\Omega}_n^{Spin}$
is  trivial, and its only
element can be represented by the sphere of the corresponding dimension
(which of course has positive Yamabe invariant).
For $n=4$, the spin cobordism group is isomorphic to ${\bf Z}$ generated
by the class of 
a $K3$ surface. Note that the Yamabe invariant of the $K3$ surface is
zero, since it admits a scalar flat metric
(the Ricci-flat metrics constructed by Yau \cite{Yau})
but no metric of positive scalar curvature (since its $\hat{A}$-genus
does not vanish). The same is true for
the connected sum of any (positive) number of copies of the
$K3$ surface. The zero
element in the group can be represented by the  4-sphere.

Now recall that D.D. Joyce  constructed examples of 8-dimensional compact
Riemannian manifolds with holonomy Spin(7) (see \cite{Joyce}). These
manifolds are then necessarily simply connected, spin,
Ricci-flat and have $\hat{A}$-genus 1 (actually D.D. Joyce shows explicitly
that his examples have $\hat{A}$-genus 1, and use this to prove that
they have holonomy Spin(7) and not a proper subgroup of it). Call $J_8$
one of these examples. Then $\alpha [J_8 ]$  is a generator of
$KO_8 (pt)$. 
Hence, multiplication by $\alpha [J_8 ]$ gives an isomorphism between
$KO_n (pt)$ and $KO_{n+8} (pt)$.  
Note also that the Yamabe invariant of $J_8$ is zero (it can not
be positive since the $\hat{A}$-genus is non-zero).

Now consider any class $[X]\in {\Omega}_n^{Spin}$ with
$n\geq 8$. Let $P$ be a closed spin manifold of dimension
$n-8$ so that $\alpha [P]$ is a generator of $KO_{n-8} (pt)$
($\alpha$ is an epimorphism in every dimension),
and let $Q=P\times J_8$. Note that it follows from Lemma 1 that
$Y(Q)\geq 0$. Note also that $\alpha [Q]$ is a generator of
$KO_n (pt)$. 
There exists then an integer $k$ so that $\alpha [X] =k\alpha [Q]$.
This means that $[X]-k[Q]$ is in the kernel of $\alpha$ and  therefore,
from Stolz's Theorem (see \cite{Stolz}), it can be represented by a closed
connected spin manifold $S$ of strictly positive Yamabe invariant.
Finally, $[X]$ can be represented by  the
manifold $\overline{X}$ obtained as the connected sum of $S$ and 
$k$ copies of $Q$. But the condition that the Yamabe invariant is 
non-negative is closed under connected sum (see \cite{Kobayashi}), and
hence $Y(\overline{X} )\geq 0$. This completes the proof of Theorem 2.

\QED

\section{Proof of Theorem 1 and final remarks}

We will now prove Theorem 1. We will use 
the following result of Gromov  and Lawson \cite[proof of Theorem B]{Gromov}:
 
\begin{theorem}: Let $N$ be a compact simply connected spin n-dimensional
manifold. Suppose that n$\geq$5 and that the manifold $X$ is spin cobordant
to $N$. Then $N$ is obtained from $X$ by  doing surgery (on $X$) on 
spheres of codimension greater than two.
\end{theorem}

Let $M$ be a compact simply connected manifold of dimension at least five.
It was proved by Gromov and Lawson in \cite{Gromov}, that if 
$M$ is not spin then it admits
a metric of positive scalar curvature. This means that the Yamabe
invariant of $M$ is strictly positive. We can therefore assume that
$M$ is spin. 

We have proved in Theorem 2 that  $M$ is spin cobordant to a spin
manifold $X$ with non-negative Yamabe invariant.  
It follows from Theorem  3 that $M$ is obtained from $X$ by
doing surgery on spheres of codimension greater than two. 
Finally, it is proved in \cite{Petean} that if $\widehat{N}$ is obtained
from the compact smooth manifold $N$ by doing surgery on 
spheres of codimension greater 
than two then $Y(\widehat{N}) \geq Y(N)$. 
Applying this  result to our situation,
we see that $Y(M)\geq Y(X) \geq 0$, and we have therefore finished the proof
of Theorem 1.

\QED
 
\noindent
{\it Remark}: It is clear from Theorem 1 that for any simply connected
compact manifold $X$ of dimension greater than four so that $\alpha [X]
\neq  0$ the Yamabe  invariant of $X$ is zero.

\vspace{.3cm}

\noindent
{\it Remark}: One of the main motivations for the study of the Yamabe invariant
is the  minimax  method to construct Einstein metrics \cite{Schoen3}. Assume
for simplicity that $Y(X)\leq 0$. Then the Yamabe invariant of $X$  is  the
supremum of the  scalar curvature  over the family of metrics on $X$ of
constant scalar curvature and unit volume. Moreover, if the supremum is
achieved  by a metric $g$ (as before), then $g$ is Einstein.  Unfortunately,
it is not always the case that the supremum is achieved. One can deduce for 
instance many examples from Theorem 1. Namely, if a simply connected manifold
$X$ (compact of dimension greater than four) has non-vanishing $\alpha$-genus,
then $Y(X)=0$. Hence  the minimax procedure would provide a Ricci-flat
metric. But, as it is well-known since the work of Lichnerowicz,
when the scalar curvature vanishes the Weitzenb\"{o}ck formula for
the Dirac operator
shows that every harmonic spinor is parallel. 
And the existence of non-trivial parallel spinors is a very restrictive
condition (we know that in our examples there are non-trivial harmonic
spinors since the $\alpha$-genus is not zero).
It implies that our simply connected manifold must be the
product of certain 8-dimensional manifolds and a Ricci-flat K\"{a}hler
manifold (see \cite[Theorem 1.2 and footnote p.54]{Hitchin}). 
This, of course, implies that for most of the examples we consider
the minimax method to construct an Einstein metric does not work.
In particular, it does not work for any of the exotic spheres $S$ 
such that $\alpha [S] \neq  0$ (see
\cite{Hitchin, Milnor}). 

\vspace{.3cm}

{\bf Acknowledgements}: The author would like to thank Claude LeBrun for 
many very
useful observations related to this work. He would also like to thank the
staff and directors of the Max Planck Institut for their hospitality
during the preparation of this work and to Vyacheslav Krushkal and
Ernesto Lupercio for many helpful discussions.

\end{document}